\documentclass[fleqn]{mat01}
\usepackage{times,mathtimy,amssymb}
\begin{document}

\setcounter{page}{355}
\firstpage{355}

\font\zz=msam10 at 10pt
\def\cd{\mbox{\zz{\char'245}}}

\def\Z{\mathbb{Z}}

\newtheorem{theore}{Theorem}
\renewcommand\thetheore{\arabic{section}.\arabic{theore}}
\newtheorem{lem}[theore]{Lemma}
\newtheorem{theor}[theore]{\bf Theorem}

\def\pott{\trivlist\item[\hskip\labelsep{{\it Proof of the Theorem.}}]}

\renewcommand{\theequation}{\thesection\arabic{equation}}

\title{$\textbf{\textit{L}}^{\mathbf 1}$-convergence of complex double Fourier
series}

\markboth{Kulwinder Kaur, S~S~Bhatia and Babu Ram}{$L^{1}$-convergence
of complex double Fourier series}

\author{KULWINDER KAUR, S~S~BHATIA and BABU RAM$^{*}$}

\address{School of Mathematics and Computer Applications,
Thapar Institute of Engineering and Technology,
Post Box No.~32, Patiala~147~004, India\\
\noindent $^{*}$Department of Mathematics,
Maharshi Dayanand University, Rohtak, India\\
\noindent Email: mathkk@hotmail.com}

\volume{113}

\mon{November}

\parts{4}

\begin{abstract}
It is proved that the complex double Fourier series of an integrable
function $f(x,y)$ with coefficients $\{ c_{jk}\}$ satisfying certain
conditions, will converge in $L^{1}$-norm. The conditions used here are
the combinations of Tauberian condition of Hardy--Karamata kind and its
limiting case. This paper extends the result of Bray [1] to complex
double Fourier series.
\end{abstract}


\maketitle

\section{Introduction}

Let $\{c_{jk}\}$ be the Fourier
coefficients of an integrable function $f(x,y)$ 
on the two-dimensional torus $T^{2} = [-\pi,\pi) \times [-\pi ,\pi)$.

We consider the double Fourier series of $f$,
\begin{equation}
\sum\limits_{j=-\infty}^{\infty}\ \sum\limits_{k=-\infty}^{\infty} c_{jk}
{\rm e}^{i(jx + ky)}
\end{equation}
with rectangular partial sums $S_{mn}(x,y)$ and Ces\`{a}ro
means $\sigma_{mn} (x,y)$ defined by
\begin{align*}
S_{mn} (f; x, y) &= \sum\limits_{\vert j\vert \leq m}\ \sum\limits_{\vert
k\vert \leq n} c_{jk} {\rm e}^{i(jx + ky)},\\[.2pc]
\sigma_{mn}(f; x,y) &= \frac{1}{(m + 1) (n +
1)}\sum\limits_{j=0}^{m}\ \sum\limits_{k=0}^{n}S_{jk}(x,y) \qquad (m,n\geq
0)
\end{align*}
and for $\lambda >1$, the de Ia Vall\'{e}e-Poussin means by\vspace{-.1pc}
\begin{equation*}
V_{mn}^{\lambda} (f, x, y) = \frac{1}{([\lambda m] -m) ([\lambda n]
-n)}\sum\limits_{j=m+1}^{[\lambda m]}\ \sum\limits_{k=n+1}^{[\lambda n]}
S_{jk}(x,y),
\end{equation*}
where $[t]$ stands for the greatest integer $\leq t$, for a positive
real number $t$.

For $f$ in $L^1(T^2), \Vert f\Vert = \int_{-\pi}^{\pi} \int_{-\pi}^{\pi}
\vert f(x,y) \vert {\rm d}x\ {\rm d}y$ stands for the $L^{1}$-norm.

The aim of this paper is to show that for a function $f$ in $L^1(T^2)$,
under certain conditions on the double Fourier series, the rectangular
partial sums $S_{mn}$ converge in $L^1$-norm to $f$, that is, $\Vert
S_{mn}-f \Vert\rightarrow 0$ as $m, n\rightarrow\infty$. Our result
belongs to the extensive study in the literature of what are called
$L^1$-convergence classes in one as well as two variables. In the single
variable case, we particularly refer to \cite{1}, \cite{2}, \cite{3}
(see \cite{2} for further references) and for double Fourier series we
are aware of \cite{4} and \cite{8}. For functions $f$ of a single
variable with the Fourier series $\sum_{n=-\infty}^{\infty}c_{n} {\rm
e}^{inx}$, a set of conditions on the coefficients $c_n$ is said to
define an $L^1$-convergence class if, under these conditions,
\begin{equation}
\Vert S_n - f\Vert = o(n)\ \hbox{if and only if}\ c_n \log n=o(n).
\end{equation}

Here $S_n$ is the $n$th partial sum of the Fourier series, $S_n =
\sum_{k=-n}^{n} c_{k} {\rm e}^{ikx}$ and $\Vert \cdot \Vert$ is the
$L^1$-norm on $[-\pi,\pi)$. The obvious modification in the definition
of $L^1$-convergence class for double Fourier series is to demand
instead of (1.2)
\begin{equation*}
\Vert S_{mn}-f\Vert \rightarrow 0\ \hbox{if and only if}\ c_{mn}\log m
\log n\rightarrow 0\ \hbox{as min}\ (m,n)\rightarrow\infty.
\end{equation*}

$\left.\right.$\vspace{-1.5pc}

Our result is weaker in not having any ``only if'' part as in the above. To
put the conditions in our result in perspective we quote the condition
for single Fourier series
\begin{equation}
\lim\limits_{\lambda \downarrow 1}\ \mathop{\overline{\lim}}\limits_{n\rightarrow \infty}\
\sum\limits_{\vert k\vert =n}^{[ \lambda n]}\vert k\vert^{p-1} \vert
\Delta c_k \vert^{p}=0,\qquad (1<p\leq 2).
\end{equation}

This condition used in \cite{1} and \cite{2} is a Tauberian condition of
the Hardy--Karamata kind \cite{3} and is a weaker version of the
conditions used by Fomin \cite{5}, Kolmogorov \cite{6}, Littlewood \cite{7} and
Teljakovskii \cite{9}. A limiting condition of (1.3) as $p\rightarrow 1$
is
\begin{equation}
\lim\limits_{\lambda \downarrow 1}\ \mathop{\overline{\lim}}\limits_{n\rightarrow
\infty} (\log n) \sum_{\vert k\vert =n}^{[\lambda n]} \vert
\Delta c_k \vert =0
\end{equation}
which can be thought of as the prototype of the conditions we use in our
result (conditions (3.1) to (3.4) below). Our conditions have points of
contact with \cite{8} and our results are simpler than those of \cite{8}
and \cite{9}. Our methods are heavily dependent on some identities in
trigonometric series involving the de la Valle\'e-Poussin means
$V^\lambda_{mn}$. The identities are developed in \S2 and the main
result is proved in \S3.

We end this section by explaining some notations.

Let
\begin{align*}
\triangle_{00} c_{jk} &= c_{jk};\\[.2pc]
\triangle_{pq} c_{jk} &= \triangle_{p-1,q} c_{jk} - \triangle_{p-1,q}
c_{j+1,k} \qquad (p\geq 1),\\[.2pc]
\triangle_{pq}c_{jk} &= \triangle_{p,q-1} c_{jk} - \triangle_{p,q-1}
c_{j,k+1} \qquad (q\geq 1).
\end{align*}
We mention that a double induction argument gives
\begin{equation*}
\triangle_{pq}c_{jk} = \sum\limits_{s=0}^{p}\ 
\sum\limits_{t=0}^{q}(-1)^{s+t}\dbinom{p}{s} \dbinom{q}{t}c_{j+s,\text{
}k+t}.
\end{equation*}

In this paper, we use the differences $\triangle_{pq}$ with $p,q\leq 2$.

For $f\in L^{1} (T^{2})$, the symbols $S_{mn}(f)$ and $S_{mn}(f,x,y)$
will have the same meaning\break as $S_{mn}$.

Similarly 
\begin{equation*}
\sigma_{mn} (f) = \sigma_{mn} (f;x,y) =\sigma_{mn}
\end{equation*}
and
\begin{equation*}
V_{mn}^{\lambda} (f) = V_{mn}^{\lambda} (f;x,y) = V_{mn}^{\lambda}.
\end{equation*}

In the sequel $\lambda_{n} = [\lambda n]$ where $n$ is a positive integer
and $\lambda >1$ is a real number.

Finally, we define the functions $E_{0+}(x) = E_{0-} (x) =\frac{1}{2},
E_{n}(x) = \sum_{k=0}^{n} {\rm e}^{ikt}, n\in \Z$ and $E_{-n}(x) =
E_{n}(-x)$. We have $\Vert E_{k}(x) \Vert \leq C\log \vert k\vert$,
where $C$ is an absolute constant.

\section{Lemmas}

In order to establish our result, we need the following lemmas:

\begin{lem}\hskip -.535pc{\rm \cite{2}}.\ \ For $m,n\geq 0$ and $\lambda >1${\rm ,} the
following representation holds{\rm :}
\begin{align*}
S_{mn}-\sigma_{mn} &=
\frac{\lambda_{m}+1}{\lambda_{m}-m}\frac{\lambda_{n}+1}{\lambda_{n}-n}
(\sigma_{\lambda_{m},\text{ }\lambda_{n}} -\sigma_{\lambda_{m},\text{ }
n}-\sigma_{{m},\text{ }\lambda_{n}} + \sigma_{mn})\\[.2pc]
&\quad\ + \frac{\lambda_{m}+1}{\lambda_{m}-m}
(\sigma_{\lambda_{m},\text{ }n}-\sigma_{mn})+
\frac{\lambda_{n}+1}{\lambda_{n}-n} (\sigma_{m,\text{ }\lambda
_{n}}-\sigma_{mn})\\[.2pc]
&\quad\ -\sum\limits_{\vert j\vert \leq m}\ \sum\limits_{\vert k\vert
=n+1}^{\lambda_{n}}\frac{\lambda_{n}+1-\vert k\vert }{\lambda
_{n}-n}c_{jk}{\rm e}^{i(jx + ky)}\\[.2pc]
&\quad\ - \sum\limits_{\vert j\vert =m+1}^{\lambda_{m}}
\ \sum\limits_{\vert k\vert \leq n} \frac{\lambda_{m}+1-\vert
j\vert}{\lambda_{m}-m} c_{jk} {\rm e}^{i(jx+ky)}\\[.2pc]
&\quad\ -\sum\limits_{\vert j\vert =m+1}^{\lambda
_{m}}\ \sum\limits_{\vert k\vert =n+1}^{\lambda_{n}}\frac{\lambda
_{m}+1-\vert j\vert }{\lambda_{m}-m}\frac{\lambda_{n}+1-\vert k\vert
}{\lambda_{n}-n}c_{jk}{\rm e}^{i(jx+ky)}.
\end{align*}
\end{lem}

\begin{lem}For $m,n\geq 0$ and $\lambda >1${\rm ,} we have the following
representation{\rm :}
\begin{align*}
V_{mn}^{\lambda}-S_{mn} &= \sum_{\vert j\vert =m+1}^{\lambda_{m}}\ 
\sum_{\vert k\vert =n+1}^{\lambda_{n}} \frac{\lambda_{m} + 1-\vert
j\vert}{\lambda_{m}-m}\frac{\lambda_{n}+1-\vert k\vert}{\lambda_{n}-n}
c_{jk} {\rm e}^{i(jx + ky)}\\[.2pc]
&\quad\ +\sum\limits_{\vert j\vert \leq m}\ \sum\limits_{\vert k\vert
=n+1}^{\lambda_{n}} \frac{\lambda_{n}+1-\vert k\vert }{\lambda
_{n}-n}c_{jk} {\rm e}^{i(jx + ky)}\\[.2pc]
&\quad\ + \sum\limits_{\vert j\vert =m+1}^{\lambda_{m}}\ \sum\limits_{\vert
k\vert \leq n} \frac{\lambda_{m}+1-\vert j\vert }{\lambda_{m}-m}c_{jk}
{\rm e}^{i(jx+ky)}.
\end{align*}
\end{lem}

\begin{proof} We have
\begin{equation*}
{V}_{mn}^{\lambda }(f,x,y) = \frac{1}{(\lambda_{m}-m) (\lambda_{n}-n)}
\sum\limits_{j=m+1}^{\lambda_{m}}\ \sum\limits_{k=n+1}^{\lambda_{n}}
S_{jk} (x,y).
\end{equation*}
Performing the double summation by parts, we have
\begin{align*}
V_{mn}^{\lambda} &= \frac{\lambda_{m}+1}{\lambda_{m}-m}
\frac{\lambda_{n}+1}{\lambda_{n}-n}\sigma_{\lambda_{m},\text{ }\lambda
_{n}}-\frac{\lambda_{m}+1}{\lambda_{m}-m}\frac{n+1}{\lambda_{n}-n}
\sigma_{\lambda_{m},\text{ }n}\\[.2pc]
&\quad\ -\frac{m+1}{\lambda_{m}-m}\frac{\lambda_{n}+1}{\lambda
_{n}-n}\sigma_{m,\text{ }\lambda_{n}}+\frac{m+1}{\lambda_{m}-m}
\frac{n+1}{\lambda_{n}-n}\sigma_{mn}\\[.2pc]
&= \frac{\lambda_{m}+1}{\lambda_{m}-m}\frac{
\lambda_{n}+1}{\lambda_{n}-n}(\sigma_{\lambda_{m},\text{ }\lambda
_{n}}-\sigma_{\lambda_{m},\text{ }n}-\sigma_{m,\text{ }\lambda
_{n}}+\sigma_{mn})\\[.2pc]
&\quad\ +\frac{\lambda_{m}+1}{\lambda_{m}-m}(\sigma _{\lambda_{m},
\text{ }n}-\sigma_{mn})+\frac{\lambda_{n}+1}{\lambda_{n}-n}
(\sigma_{\lambda_{m},\text{ }n}-\sigma_{mn}) +\sigma_{mn}
\end{align*}
The use of Lemma~2.1, gives
\begin{align*}
V_{mn}^{\lambda} -S_{mn} &= \underset{\vert j\vert
=m+1}{\overset{\lambda_{m}}{\sum }}\ \overset{\lambda_{n}}{
\underset{\vert k\vert =n+1}{\sum }}\frac{\lambda_{m}+1-\vert j\vert
}{\lambda_{m}-m}\frac{\lambda_{n}+1-\vert k\vert }{\lambda
_{n}-n}\mathit{c}_{jk} {\rm e}^{i(jx+ky)}\\[.2pc]
&\quad\ +\sum\limits_{\vert j\vert \leq m}\
\sum\limits_{\vert k\vert =n+1}^{\lambda_{n}} \frac{\lambda
_{n}+1-\vert k\vert }{\lambda_{n}-n}c_{jk}{\rm e}^{i(jx+ky)}\\[.2pc]
&\quad\ - \sum\limits_{\vert j\vert =m+1}^{\lambda_{m}}\ \sum\limits_{\vert k\vert \leq
n} \frac{\lambda _{m}+1-\vert j\vert
}{\lambda_{m}-m}c_{jk}{\rm e}^{i(jx+ky)}.
\end{align*}
\end{proof}

\begin{lem} For $m,n\geq 0$ and $\lambda >1${\rm ,} we have
\begin{align*}
&\sum\limits_{j\leq m}\ \sum\limits_{\vert k\vert
=n+1}^{\lambda_{n}}\frac{\lambda_{n}+1-\vert k\vert }{\lambda_{n}-n}c_{jk}{\rm
e}^{i(jx+ky)}\\[.2pc]
&= \sum\limits_{j=0\pm }^{m-1}\ \sum\limits_{\vert k\vert
=n}^{\lambda_{n}-1}\frac{\lambda_{n}-\vert
k\vert}{\lambda_{n}-n}\Delta_{11}c_{jk}E_{j}(x) E_{k}(y)\\[.2pc]
&\quad\ +\frac{1}{\lambda_{n}-n}\overset{m-1}{\sum\limits_{\vert j\vert
=\pm 0}}\ \sum\limits_{\vert k\vert
=n+1}^{\lambda_{n}}\Delta_{10}c_{jk}E_{j}(x) E_{k}(y)\\[.2pc]
&\quad\ -\sum\limits_{\vert j\vert =0\pm }^{\lambda
_{n}}\ \sum\limits_{k=n}\Delta_{10}c_{jk}E_{j}(x) E_{k}(y)\\[.2pc]
&\quad\ +\sum\limits_{j=m}\ \sum\limits_{\vert k\vert =n}^{\lambda
_{n}-1}\frac{\lambda_{n}-\vert k\vert }{\lambda_{n}-n}\Delta
_{01}c_{jk}E_{j}(x) E_{k}(y)\\[.2pc]
&\quad\ +\frac{1}{\lambda_{n}-n}\sum\limits_{\vert j\vert =m}
\ \sum\limits_{\vert k\vert =n+1}^{\lambda_{n}} c_{jk}E_{j}(x)
E_{k}(y)\\[.2pc]
&\quad\ -\sum\limits_{\vert j\vert =m}\ \sum\limits_{\vert
k\vert =n}c_{jk}E_{j}(x) E_{k}(y).
\end{align*}
\end{lem}

\begin{proof}By summation by parts,
\begin{align*}
&\sum\limits_{k=n+1}^{\lambda_{n}}\frac{\lambda_{n}+1-\vert
k\vert }{\lambda_{n}-n}c_{jk}{\rm e}^{iky}\\[.2pc]
&= \overset{\lambda_{n-1}}{\sum\limits_{k=n}}\Delta_{01}
\left[\frac{\lambda_{n}+1-k}{\lambda_{n}-n}c_{jk} \right] E_{k}(y) +
\frac{1}{\lambda_{n}-n}c_{j,\text{ }\lambda_{n}}E_{\lambda_{n}}
(y)\\[.2pc]
&\quad\ -\frac{\lambda_{n}-n+1}{\lambda_{n}-n}c_{j,\text{ }n}E_{n}(y)\\[.2pc]
&=\overset{\lambda_{n-1}}{\sum\limits_{k=n}}\frac{\lambda_{n}-k}{\lambda
_{n}-n}\Delta_{01}c_{jk}E_{k}(y) +\frac{1}{\lambda_{n}-n}
\underset{k=n+1}{\overset{\lambda_{n}}{\sum }}c_{jk}E_{k}(y)
-c_{j,\text{ }n}E_{n}(y).
\end{align*}
Similarly
\begin{align*}
&\underset{k=n+1}{\overset{\lambda_{n}}{\sum }}\frac{\lambda_{n}+1-k}{
\lambda_{n}-n}c_{j,-k}{\rm e}^{-iky}\\[.2pc]
&=\underset{k=n}{\overset{\lambda_{n-1}}{\sum }}\frac{\lambda_{n}-k}{
\lambda_{n}-n}\Delta_{01}c_{j,\text{ }-k}E_{-k}(y) +\frac{1}{
\lambda_{n}-n}\underset{k=n+1}{\overset{\lambda_{n}}{\sum }}c_{j,\text{ }
-k}E_{-k}(y)\\[.2pc]
&\quad\ -c_{j,\text{ }-n}E_{-n}(y).
\end{align*}
Combining the above results
\begin{align*}
&\underset{\vert k\vert =n+1}{\overset{\lambda_{n}}{\sum }
}\frac{\lambda_{n}+1-\vert k\vert }{\lambda_{n}-n}c_{jk}{\rm e}^{iky}\\[.2pc]
&=\underset{\vert k\vert =n}{\overset{\lambda_{n}-1}{\sum
}}\frac{\lambda_{n}-\vert k\vert }{\lambda_{n}-n}\Delta
_{01}c_{jk}E_{k}(y) +\frac{1}{\lambda_{n}-n}\underset{ \vert k\vert
=n+1}{\overset{\lambda_{n}}{\sum }} c_{jk}E_{k}(y)\\[.2pc]
&\quad\ -\sum\limits_{k=\vert n\vert }c_{jk}E_{k}(y).
\end{align*}
Another summation by parts with respect to $j$ gives the required
result.
\end{proof}

\begin{lem}For $m,n\geq 0$ and $\lambda >1${\rm ,} we have
\begin{align*}
&\underset{\vert j\vert =m+1}{\overset{\lambda_{m}}{\sum }
}\ \underset{\vert k\vert =n+1}{\overset{\lambda_{n}}{\sum }}
\frac{\lambda_{m}+1-\vert j\vert }{\lambda_{m}-m}\frac{
\lambda_{n}+1-\vert k\vert }{\lambda_{n}-n}c_{jk}{\rm e}^{i(jx+ky)}\\[.2pc]
&= \underset{\vert j\vert =m}{\overset{\lambda_{m}-1} {\sum
}}\ \underset{\vert k\vert =n}{\overset{\lambda_{n}-1}{\sum }}
\frac{\lambda_{m}-\vert j\vert }{\lambda_{m}-m}\frac{\lambda _{n}-\vert
k\vert }{\lambda_{n}-n}\Delta_{11}c_{jk}E_{k}(x) E_{k}(y)\\[.2pc]
&\quad\ +\frac{1}{\lambda_{n}-n}\underset{\vert j\vert =m}{
\overset{\lambda_{m-1}}{\sum }}\ \underset{\vert k\vert =n+1}{
\overset{\lambda_{n}}{\sum }}\frac{\lambda_{m}-\vert j\vert}{
\lambda_{m}-m}\Delta_{10}c_{jk}E_{j}(x) E_{k}(y)\\[.2pc]
&\quad\ +\frac{1}{\lambda_{m}-m}\underset{\vert j\vert =m+1}{
\overset{\lambda_{m}}{\sum }}\ \underset{\vert k\vert =n}{\overset{
\lambda_{n-1}}{\sum }}\frac{\lambda_{n}-\vert k\vert }{\lambda
_{m}-m}\Delta_{01}c_{jk}E_{j}(x) E_{k}(y)
\end{align*}
\begin{align*}
&\quad\ -\underset{\vert j\vert =m}{\overset{\lambda_{m}-1}{\sum
}}\ \sum\limits_{\vert k\vert =n}\frac{\lambda_{m}-\vert j\vert
}{\lambda_{m}-m}\Delta_{10}c_{jk}E_{j}(x) E_{k}(y)\\[.2pc]
&\quad\ -\sum\limits_{j=m}\ \sum\limits_{k=n}^{\lambda_{n}-1}\frac{\lambda
_{m}-\vert k\vert }{\lambda_{n}-n}\Delta_{01}c_{jk}E_{j}(x) E_{k}(y)\\[.2pc]
&\quad\ -\frac{1}{\lambda_{m}-m}\underset{\vert j\vert =m+1}{
\overset{\lambda_{m}}{\sum }}\ \underset{\vert k\vert =n}{\overset{ }{\sum
}}c_{jk}E_{j}(x) E_{k}(y)\\[.2pc]
&\quad\ -\frac{1}{\lambda_{n}-n}\underset{\vert j\vert =m}{
\overset{}{\sum }}\ \underset{\vert k\vert =n+1}{\overset{\lambda
_{n}}{\sum }}c_{jk}E_{j}(x) E_{k}(y)\\[.2pc]
&\quad\ +\frac{1}{\lambda_{n}-n}\frac{1}{\lambda_{m}-m}\underset{ \vert
j\vert =m+1}{\overset{\lambda_{m}}{\sum }}\ \underset{ \vert k\vert
=n+1}{\overset{\lambda_{n}}{\sum }} c_{jk}E_{j}(x) E_{k}(y)\\[.2pc]
&\quad\ -\underset{\vert j\vert =m}{\overset{}{\sum }}\ \underset{ \vert
k\vert =n}{\overset{}{\sum }}c_{jk}E_{j}(x) E_{k}(y).
\end{align*}
\end{lem}

Proof of this lemma follows on similar lines as in Lemma~2.3.

\section{Main result}

The main result of this paper is the following theorem:

\setcounter{theore}{0}
\begin{theor}[\!] Let $f\in L^{1}(T^{2})${\rm ,} and $\{c_{jk}\}$ be
its Fourier coefficients. If
\setcounter{equation}{0}
\begin{align}
&\underset{\vert j\vert =0\pm}{\overset{\infty}{\sum}}(\log \vert
j\vert) (\log \vert k\vert) \vert \Delta_{10}c_{jk}\vert \rightarrow 0,\
\ {\rm as}\quad \vert k\vert \rightarrow \infty,\\[.2pc]
&\underset{\vert k\vert =0\pm }{\overset{\infty}{\sum}}(\log \vert
j\vert) (\log \vert k\vert) \vert \Delta_{01}c_{jk}\vert \rightarrow 0,\
\ {\rm as}\quad \vert j\vert \rightarrow \infty,\\[.2pc]
&\underset{\lambda \downarrow 1}{\lim }\ \mathop{\overline{\lim}}
\limits_{n\rightarrow \infty}\ \underset{\vert j\vert =0\pm
}{\overset{\infty }{\sum }}\ \underset{\vert k\vert =n}{\text{
}\overset{\lambda_{n}}{\sum }}(\log \vert j\vert) (\log \vert k\vert
) \vert \Delta _{11}c_{jk}\vert =0,\\[.2pc]
&\underset{\lambda \downarrow 1}{\lim }\textit{\ }
\mathop{\overline{\lim}}\limits_{m\rightarrow \infty}\ \underset{\vert
k\vert =0\pm }{\overset{\infty }{\sum }}\ \underset{\vert j\vert
=m}{\text{ }\overset{\lambda_{m}}{\sum }}(\log \vert j\vert) (\log
\vert k\vert) \vert \Delta_{11}c_{jk}\vert =0.
\end{align}
Then $\Vert S_{mn}(f) -f\Vert =o(1)$ as $\min (m,n) \rightarrow
\infty.$
\end{theor}

\begin{pott}{\rm Since $f\in L^{1}(T^{2})$ therefore $\Vert
\sigma_{mn}(f) -f\Vert =o(1)$ as $\min (m,n)$ $\rightarrow \infty$ (see,
e.g.~\cite{10}, vol.~2, p.~309). It follows that $\Vert
V_{mn}^{\lambda}-f\Vert =o(1)$ as $(m,n) \rightarrow \infty$.

Consequently it is sufficient to prove that
\begin{equation*}
\Vert V_{mn}^{\lambda }-S_{mn}\Vert =o(1) \ \ {\rm as}\ \ \min (m,n)
\rightarrow \infty.
\end{equation*}}
\end{pott}

Combining the results of Lemmas~2.1--2.4, we have
\begin{align*}
V_{mn}^{\lambda}-S_{mn} &= R_{1}^{\lambda }(m,n;x,y) +R_{2}^{\lambda} 
(m,n;x,y) -R_{3}^{\lambda }(m,n;x,y)\\[.2pc]
&\quad\ -R_{4}^{\lambda }(m,n;x,y) + R_{5}^{\lambda }(m,n;x,y)
-R_{0}^{\lambda }(m,n;x,y),\\[.2pc]
R_{1}^{\lambda }(m,n;x,y) &=\sum\limits_{\vert j\vert
=m}^{\lambda_{m}-1}\ \sum\limits_{\vert k\vert =n}^{\lambda_{n}-1}
\frac{\lambda_{m}-\vert j\vert }{\lambda_{m}-m}\frac{\lambda _{n}-\vert
k\vert }{\lambda_{n}-n}\Delta_{11}c_{jk}E_{j}(x) E_{k}(y)\\[.2pc]
&\quad\ +\sum\limits_{\vert j\vert =0\pm }^{m-1}\ \sum\limits_{\vert
k\vert =n}^{\lambda_{n}-1}\frac{ \lambda_{n}-\vert k\vert
}{\lambda_{n}-n}\Delta _{11}c_{jk}E_{j}(x) E_{k}(y)\\[.2pc]
&\quad\ +\sum\limits_{\vert j\vert =m}^{\lambda_{m}-1}\ \sum\limits_{\vert
k\vert =0\pm }^{n-1}\frac{ \lambda_{m}-\vert j\vert
}{\lambda_{m}-m}\Delta _{11}c_{jk}E_{j}(x) E_{k}(y),\\[.2pc]
R_{2}^{\lambda}(m,n;x,y) &=\frac{1}{\lambda_{n}-n} \sum\limits_{\vert
j\vert =0\pm }^{m-1}\ \sum\limits_{\vert k\vert
=n+1}^{\lambda_{n}-1}\Delta_{10}c_{jk}E_{j}(x) E_{k}(y)\\[.2pc]
&\quad\ +\frac{1}{\lambda_{m}-m}\sum\limits_{\vert j\vert
=m+1}^{\lambda_{m}}\ \sum\limits_{\vert k\vert =0\pm
}^{n-1}\Delta_{01}c_{jk}E_{j}(x) E_{k}(y)\\[.2pc]
&\quad\ +\frac{1}{\lambda_{n}-n}\sum\limits_{\vert j\vert
=m+1}^{\lambda_{m}}\ \sum\limits_{\vert k\vert
=n}^{\lambda_{n}-1}\frac{\lambda_{m}-\vert j\vert }{\lambda
_{m}-m}\Delta_{10}c_{jk}E_{j}(x) E_{k}(y)\\[.2pc]
&\quad\ +\frac{1}{\lambda_{m}-m}\sum\limits_{\vert j\vert
=m+1}^{\lambda_{m}}\ \sum\limits_{\vert k\vert
=n}^{\lambda_{n}-1}\frac{\lambda_{n}-\vert j\vert }{\lambda
_{n}-n}\Delta_{01}c_{jk}E_{j}(x) E_{k}(y),\\[.2pc]
R_{3}^{\lambda }(m,n;x,y) &=\sum\limits_{\vert k\vert
=n}\ \sum\limits_{\vert j\vert =0\pm }^{m-1}\Delta _{10}c_{jk}E_{j}(x)
E_{k}(y),\\[.2pc]
R_{4}^{\lambda }(m,n;x,y) &=\sum\limits_{\vert j\vert
=m}\ \sum\limits_{\vert k\vert =0\pm }^{n-1}\Delta _{01}c_{jk}E_{j}(x)
E_{k}(y),\\[.2pc]
R_{5}^{\lambda }(m,n;x,y) &=\frac{1}{\lambda_{m}-m}\frac{1}{
\lambda_{n}-n}\sum\limits_{\vert j\vert =m+1}^{\lambda
_{m}}\ \sum\limits_{\vert k\vert =n+1}^{\lambda _{n}}c_{jk}E_{j}(x)
E_{k}(y),\\[.2pc]
R_{0}^{\lambda }(m,n;x,y) &=\sum\limits_{\vert j\vert
=m}\ \sum\limits_{\vert k\vert =n}c_{jk}E_{j}(x) E_{k}(y).
\end{align*}
It follows from (1.6),
\begin{align*}
\Vert R_{1}^{\lambda}(m,n;x,y) \Vert_{1} &\leq
C \left\{\sum\limits_{\vert j\vert =0\pm }^{\infty
}\ \sum\limits_{\vert k\vert =n}^{\lambda_{n}}(\log
\vert j\vert) (\log \vert k\vert)
\vert \Delta_{11}c_{jk}\vert \right.\\[.2pc]
&\quad\ \left. +\sum\limits_{\vert j\vert
=0\pm }^{\infty }\ \sum\limits_{\vert k\vert =n}^{\lambda
_{n}}(\log \vert j\vert) (\log \vert
k\vert) \vert \Delta_{11}c_{jk}\vert \right\}.
\end{align*}
Thus, by (3.3) and (3.4), we get
\begin{equation*}
\underset{\lambda \downarrow 1}{\lim } {\text{ 
}} \mathop{\overline{\lim}}\limits_{m,n\rightarrow \infty} \Vert R_{1}^{\lambda }(
m,n;x,y) \Vert_{1}=0.
\end{equation*}
$R_{2}^{\lambda}(m,n;x,y) = R_{21}^{\lambda}(m,n;x,y) +
R_{22}^{\lambda}(m,n;x,y)$, where
\begin{align*}
R_{21}^{\lambda}(m,n;x,y) &=\frac{1}{\lambda_{n}-n} \sum\limits_{\vert
j\vert =0\pm }^{m-1}\ \sum\limits_{\vert k\vert
=n+1}^{\lambda_{n}}\Delta_{10}c_{jk}E_{j}(x) E_{k}(y)\\[.2pc]
&\quad\ +\frac{1}{\lambda_{n}-n}
\sum\limits_{j=m}^{\lambda_{m}-1}\ \sum\limits_{\vert k\vert
=n+1}^{\lambda_{n}}\frac{\lambda_{m}-\vert j\vert }{\lambda
_{m}-m}\Delta_{10}c_{jk}E_{j}(x) E_{k}(y),\\[.2pc]
R_{22}^{\lambda}(m,n;x,y) &=\frac{1}{\lambda_{m}-m} \sum\limits_{\vert
j\vert =m+1}^{\lambda _{m}}\ \sum\limits_{\vert k\vert =\pm 0}^{n-1}\Delta
_{01}c_{jk}E_{j}(x) E_{k}(y)\\[.2pc]
&\quad\ +\frac{1}{\lambda_{m}-m}\sum\limits_{ \vert j\vert
=m+1}^{\lambda_{m}}\ \sum\limits_{\vert k\vert
=n}^{\lambda_{n}-1}\frac{\lambda_{n}-\vert k\vert }{\lambda
_{n}-n}\Delta_{01}c_{jk}E_{j}(x) E_{k}(y),\\[.2pc]
R_{21}^{\lambda}(m,n;x,y) &=\frac{1}{\lambda_{n}-n} \sum\limits_{\vert
j\vert =0}^{m-1}\ \sum\limits_{\vert k\vert
=n+1}^{\lambda_{n}}\Delta_{10}c_{jk}E_{j}(x) E_{k}(y)\\[.2pc]
&\quad\ +\frac{1}{\lambda_{n}-n}\sum\limits_{ \vert j\vert
=m}^{\lambda_{m-1}}\ \sum\limits_{\vert k\vert
=n+1}^{\lambda_{n}}\frac{\lambda_{m}-\vert j\vert }{\lambda
_{m}-m}\Delta_{10}c_{jk}E_{j}(x) E_{k}(y)\\[.2pc]
&=\frac{1}{(\lambda_{m}-m) (\lambda _{n}-n) }\\[.2pc]
&\quad\ \times \sum\limits_{\vert i\vert =m}^{\lambda _{m-1}}\
\sum\limits_{\vert k\vert =n+1}^{\lambda_{n}} \left(\sum\limits_{\vert
j\vert =0}^{i}\Delta_{10}c_{jk}E_{j}(x) E_{k}(y)\right).
\end{align*}
By (1.6) and (3.1), for $\lambda >1$, we conclude that
\begin{equation*}
\underset{\lambda \downarrow 1}{\lim }\ 
\mathop{\overline{\lim}}\limits_{m,n\rightarrow \infty} \Vert R_{21}^{\lambda}
(m,n;x,y) \Vert_{1}=0.
\end{equation*}
Similarly, by (1.6) and (3.2), for $\lambda >1$, we have
\begin{equation*}
\underset{\lambda \downarrow 1}{\lim }\
\mathop{\overline{\lim}}\limits_{m,n\rightarrow \infty} \Vert
R_{22}^{\lambda} (m,n;x,y) \Vert_{1}=0.
\end{equation*}
By (3.1) and (3.2), we get
\begin{align*}
\underset{\lambda \downarrow 1}{\lim } 
\mathop{\overline{\lim}}\limits_{m,n\rightarrow \infty} \Vert
R_{3}^{\lambda }(m,n;x,y) \Vert_{1} &= 0,\\[.2pc]
\underset{\lambda \downarrow 1}{\lim }\ 
\mathop{\overline{\lim}}\limits_{m,n\rightarrow \infty} \Vert
R_{4}^{\lambda}(m,n;x,y) \Vert_{1} &= 0.
\end{align*}
Taking into account $\vert c_{jk}\vert \leq \sum_{\vert u\vert
=j}^{\lambda_{n}}\vert \Delta _{10}c_{uk}\vert$, and by (3.1), we find
that $c_{jk}\log \vert j\vert \log \vert k\vert =o(1)$ as $\min (\vert
j\vert,\vert k\vert) \rightarrow \infty$.

Combining this with (1.6), we conclude that
\begin{equation*}
\vert R_{5}^{\lambda }(m,n;x,y) \vert \leq \underset{
\underset{n<v\leq \lambda_{n}}{m<u\leq \lambda_{m}}}{\mathrm{max}}
\Vert R_{0}^{\lambda }(u,v;x,y) \Vert \rightarrow
\infty\ \ {\rm as}\ \ \min (m,n) \rightarrow \infty.
\end{equation*}

$\left.\right.$\vspace{-1.6pc}

Combining all what we have done so far, we conclude the desired result.

\section*{Acknowledgement}

The authors are thankful to the referee for suggestions which have
improved the presentation of the paper and for calling their attention to
the paper of Moricz \cite{8}.

\end{document}